\documentstyle[12pt]{article}
\newcommand{\be}{\begin{equation}}
\newcommand{\ee}{\end{equation}}
\newcommand{\bqn}{\begin{eqnarray}}

\newcommand{\eqn}{\end{eqnarray}}

\newcommand{\bd}{\begin{description}}
\newcommand{\ed}{\end{description}}

\newtheorem{stat}{}[section]

\def\bs{\begin{stat}}
\def\es{\end{stat}}
\def\ben{\begin{enumerate}}
\def\een{\end{enumerate}}

\def\bp{\noindent{\bf Proof}  \ \ \ }
\def\ep{\hfill $\Box$}

\makeatletter
\@addtoreset{equation}{section}

\makeatother

\begin{document}

\begin{center}
{\large {\bf SIMPLE AND DIRECT PROOF }}
\\[2ex]
{\large {\bf OF MACLANE's GRAPH PLANARITY CRITERION}}
% \\[2ex]
% {\large {\bf for Graphs }}
\\[4ex]
{\large {\bf Alexander Kelmans}}
\\[2ex]
{\bf Rutgers University, New Brunswick, New Jersey}
%\\[0.5ex]
%{\bf and}
\\[0.5ex]
{\bf University of Puerto Rico, San Juan, Puerto Rico}
\\[2ex]
\end{center}

% \begin{abstract}
\abstract
{We give a simple proof of 
%(a natural refinement of) 
MacLane's algebraic planarity criterion for graphs. This proof does not use any other known planarity criteria.
\\[1ex]
}
% \end{abstract}

%\rrr{??-2005}

\section{Introduction}

\indent

We consider undirected graphs with no loops (parallel edges are possible). All notions on graphs, that are  not defined here, can be found in \cite{D,V}.
\\

There are various graph planarity criteria.
Here are some of them.

\bs{\em (Kuratowski \cite{KK})}
\label{Kuratowski}
A graph is non-planar if and only if it contains a subdivision of $K_5$ or $K_{3,3}$.
\es

\bs {\em (Whitney \cite{W2})}
\label{Whitney} A graph is planar if and only if it has a matroid dual graph.
\es

\bs {\em (MacLane \cite{ML})}
\label{MacLane} A graph is planar if and only if its cycle space has a  2--basis $($i.e. a basis that consists of some cycles of the graph and such that  every edge of the graph belongs to at most two cycles from the basis$)$.
\es

A cycle $C$ in a connected graph $G$ is called {\em separating} if $G / C$ has more blocks then $G$, and
{\em non-separating}, otherwise.

\bs {\em (Kelmans \cite{K1,K2})}
\label{Kelmans} A 3--connected graph is planar if and only if each edge of the graph belongs to exactly two non-separating cycles of the graph.
\es

There are several fairly simple proofs of {\bf \ref{Kuratowski}}
(e.g. \cite{K2,M,Thms}).
Theorems {\bf \ref{Whitney}} and {\bf \ref{MacLane}} follow 
% easily 
from {\bf \ref{Kuratowski}} because $K_5$ and $K_{3,3}$ have  no matroid dual graph and have 
no 2--basis, respectively (e.g. \cite{D,ML,W2}).
In  \cite{K1,K2} we gave a simple proof of {\bf \ref{Kelmans}}
that does not use any other known planarity criteria.
We also gave a simple proof of  {\bf \ref{Whitney}} using 
 {\bf \ref{Kelmans}}.
Moreover, we  showed that 

\bs 
\label{3edge}{\em \cite{K1}}
A 3--connected graph has an edge
belonging to at least three non-separating cycles if and only if it has a subdivision of $K_5$ or $K_{3,3}$.
\es

This fact implies that  {\bf \ref{Kuratowski}} follows from  
{\bf \ref{Kelmans}} and vise versa and that
{\bf \ref{MacLane}} follows from {\bf \ref{Kelmans}}.

The following theorem, due to W. Tutte \cite{T}  and, independently, A. Kelmans \cite{K1,K2},
is an important result in the study of the graph cycle 
spaces. 
%of graphs.
\bs
\label{NCgenerator}
The set of non-separating circuits of a 3--connected 
graph generates the cycle space of the graph.
\es

In \cite{K1} we noted that {\bf \ref{Kelmans}} follows from 
{\bf \ref{MacLane}} and  {\bf \ref{NCgenerator}}.
\\

In this paper we give a simple proof of (a natural refinement of) MacLane's graph planarity criterion {\bf \ref{MacLane}}. This proof does not use any other known planarity criteria.
\\

More information on this topic  (in particular, some strengthenings of {\bf \ref{Kuratowski}}, {\bf \ref{Whitney}}, and {\bf \ref{Kelmans}})
can be found in the expository paper \cite{K3} and in \cite{K4}.
\\

The results of this paper were presented at the Moscow Discrete Mathematics Seminar in 1977 (see also \cite{Krrr3-06}).
% In this paper we give a simple proof of MacLane's planarity 
% criterion {\bf \ref{MacLane}} that does not use any other 
% known planarity criteria.
% \\

\section{Main notions and notation}

\indent

Let $G$ be a graph, $V(G)$  and $E = E(G)$ the sets of vertices and edges of $G$, respectively. Let $e(G) = |E(G)|$.
If $C$ is a cycle of $G$, then $E(C)$ is called 
a {\em circuit} of $G$.
If $X, Y \subseteq E$, then let $X + Y$ denote 
the symmetric difference of $X$ and $Y$, 
i.e. $X + Y = (X \cup Y) \setminus (X \cap Y)$.
Then $2^E$ forms a vector space over $GF(2)$.
Let ${\cal C}(G)$ denote the set of 
% the edge sets of cycles 
circuits of $G$, and so ${\cal C}(G) \subseteq 2^E$. 
Let ${\cal CS}(G)$ denote the subspace of $2^E$ 
generated by ${\cal C}(G)$.
This subspace is called the {\em cycle space of} $G$.
Obviously $X \in {\cal CS}(G)$ if  and only if every 
vertex $v$ in the subgraph of $G$, induced by $X$, 
has even degree.
In particular, $\emptyset \in {\cal CS}(G)$.
%
% If  $A$ and $B$ are  subgraphs of $G$,
% we write, for simplicity,
% $G / A$ instead of $G/ E(A)$,  
% $A + B$ instead of $E(A) + E(B)$, and 
% $A \in {\cal F}$ instead of $E(A) \in {\cal F}$ for 
% ${\cal F} \subseteq 2^E$,
%
% A cycle $C$ (the corresponding circuit $E(C)$) in a 
% connected graph $G$ is called {\em separating}
% if $G /C$ has more blocks then $G$, and
% {\em non-separating}, otherwise.
% Let ${\cal NC}(G)$ denote the set of non-separating 
% circuits of $G$, and so 
% ${\cal NC}(G) \subseteq {\cal C}(G)$.
%
A basis $B$ of ${\cal CS}(G)$ is called {\em simple} 
if every edge of $G$ belongs to at most two members (edge sets) of $B$.

If ${\cal F} \subseteq 2^E$ and $H$ is a subgraph of $G$, we write $H \in {\cal F}$ and ${\cal F} \setminus \{H\}$ instead of $E(H) \in {\cal F}$ and ${\cal F} \setminus \{E(H)\}$, respectively.

If $X \subseteq E(G)$, then let $\dot{X}$ denote the subgraph of $G$ induced by $X$.

If $H$ is a plane 2--connected graph, 
then let ${\cal F}(H)$ be
% denote 
the set of facial circuits of $H$.

A path $P$  with end-vertices $x$ and $y$ is called 
a {\em path-chord of} a cycle $C$ (and of the corresponding circuit $E(C)$) in $G$ if  $V(C) \cap V(P) = \{x,y\}$, and 
$E(C) \cap E(P) = \emptyset $.

A {\em thread} in  $G$ is a path $T$ in $G$ such that 
the degree of every inner vertex of $T$ 
is equal to  two and the degree of every end-vertex of $T$ 
is not equal to two in $G$. 
Obviously, if $C$ is a cycle of $G$ and 
$E(C) \cap E(T) \ne \emptyset $, then $T \subseteq C$.
If $T$ is a thread in $G$, we write $G - (T)$ instead of 
$G - (T - End(T))$.

\section{Proof of MacLane's planarity criterion}

\indent

% A subset $B$ of the cycle space of $G$ is called {\em simple}
% if every edge of $G$ belongs to at most two members (edge % sets) of $B$.

% If $H$ is a plane 2--connected graph 
% then let ${\cal F}(H)$ denote the set 
% the edge sets of facial cycles of $H$.

% \bs {\em (MacLane \cite{ML})}
% \label{mclane}
% The following are equivalent:
% \\[1ex]
% $(a1)$ $G$ is planar and
% \\[1ex]
% $a2)$ $G$ has a simple cycle basis.
% \es
It is easy to see the following.
\bs
\label{EarAssembly}
Let $G$ be a  2--connected graph and $G$ not a cycle.
Then $G$ has a thread $T$ such that $G - (T)$ is 
a 2--connected graph.
\es

Obviously
\bs
\label{Gplane}
Let $G$ be a 2--connected planar graph,
$G_{\epsilon }$ be an embedding of $G$ into the plane,
and $F$ a facial circuit
%cycle 
of $G_{\epsilon }$.
Then ${\cal F}(G) \setminus \{F\}$ is a simple basis of 
${\cal CS}(G)$.
\es

\bs
\label{2con}
Let $G$ be a  2--connected graph and $G$ not a cycle.
If $B$ is a simple basis of ${\cal CS}(G)$, 
then $G$ is planar and there is an embedding $G_{\epsilon }$ 
of $G$ such that
$B = {\cal F}(G) \setminus \{F\}$ for some $F \in {\cal F}(G)$.
\es

\bp 
% Let $E(G) = E$.
We prove our claim by induction on $e(G)$.
If $e(G) = 3$, then our claim is obviously true.
So let $e(G) \ge 4$. By {\bf \ref{EarAssembly}}, there is a thread $T$ of $G$ such that $G' = G - (T)$ is 2--connected. 
Since $G$ is 2--connected, $T$ belongs to a cycle of $G$.
Therefore $E(T)$ belongs to at least one member of $B$.
Since $B$ is a simple basis of ${\cal CS}(G)$, 
$E(T)$ belongs to at most two members of $B$.

If $E(T)$ belongs to exactly one member of $B$, say $C$, 
then let $B': = B \setminus \{C\}$.
If $E(T)$ belongs to (exactly) two members, say $S$ and $Z$, of $B$, 
% (i.e. $E(T)\subseteq S\cap Z$), 
then let $B': = B \setminus \{S,Z\} \cup \{S+Z\}$.
Then $B'$  is a simple basis of $G'$.

% If $T$ belongs to exactly one member of $B$, say $C$, 
% then $ B': = B - \{C\}$ is a simple basis of $G'$.
% If $T$ belongs to (exactly) two members, say $S$ and $Z$, % of $B$ 
% (i.e. $E(T)\subseteq S\cap Z$), then $B': = B - \{S,Z\} \cup 
% \{S+Z\}$ is a simple basis of $G'$.

By the induction hypothesis, $G'$ is planar and there is an embedding $G'_{\alpha }$ of $G'$ such that 
$B'  = {\cal F}(G'_{\alpha }) \setminus \{D\}$ for some  
$D \in {\cal F}(G'_{\alpha })$, and so every member of $B'$ 
is a facial circuit of $G'_{\alpha }$ and every edge in 
$E(G) \setminus D$ belongs to exactly 
two facial circuits of $G'_{\alpha }$ that are members of $B'$.
 
Suppose that $B' = B \setminus \{C\}$. 
Since $B$ is a simple basis of $G$ and $B'$ is a subset of $B$, clearly $C \setminus E(T)$ is a subset of $D$. 
Since $T$ is a thread in $G$ and $C$ is an element of the cycle space of $G$, clearly $\dot{C} - (T)$ is a path, and so $\dot{C}$ is a cycle in $G$ and $T$ is a path-chord of cycle $\dot{D}$. 
Now since $D$ is a facial circuit of $G'_{\alpha }$, we can embed $T$ in the face, bounded by $\dot{D}$, to obtain from $G'_{\alpha }$ an embedding $G_{\epsilon}$ of $G$, and so $G$ is planar and 
$B = {\cal F}(G_{\epsilon }) \setminus \{C'\}$, where $\dot{C}'$ is the cycle in $\dot{D} \cup T$ containing $T$ and distinct from $\dot{C}$.

Now suppose that $B' = B \setminus \{S,Z\} \cup \{S+Z\}$. 
Then $S+Z$ is a facial circuit of $G'_{\alpha }$ which is a member of $B'$. 
We know that $E(T)\subseteq S \cap Z$.
Suppose that there is $e \in (S \cap Z) \setminus E(T)$.
% and so $e \not \in S+Z$. 
Then $e \in E(G') \setminus (S + Z)$, and so $e$ belongs to a member, say $R$, of $B'$.
Therefore $e$ belongs to three members of $B$, namely, $R$, $S$, and $Z$, 
and so $B$ is not a simple basis of ${\cal CS}(G)$, 
a contradiction.
Thus $S \cap Z = E(T)$, and so 
$S+Z = (S \cup Z )\setminus E(T)$ and $T$ is a path-chord 
of facial circuit $S+Z$ of $G'_{\alpha }$.
Then we can embed $T$ in the face, bounded by $S+Z$, 
to obtain from $G'_{\alpha }$ an embedding $G_{\epsilon}$ of $G$, and so $G$ is planar and 
${\cal F}(G_{\epsilon }) = 
{\cal F}(G'_{\alpha }) \setminus  \{S+Z\}) \cup \{S,Z\}$.

If $D = S+Z$ then ${\cal F}(G_{\epsilon }) = B$, and so
the sum of members of $B$ is equal to $\emptyset $.
Therefore $B$ is not a basis of ${\cal CS}(G)$, a contradiction.
Thus $D \ne S+Z$, and so 
$B = {\cal  F}(G_{\epsilon }) \setminus \{D\}$.
%
% We know that $G$ is 2--connected, $T$ is a path-chord of 
% a facial cycle $S+Z$ of $G'_{\alpha }$,  and $D$ is a facial 
% cycle of $G'_{\alpha }$. Therefore $T$ has at most one 
% vertex (namely, an end-vertex) in $D$. Since $D$ is a facial % cycle of $G'_{\alpha }$, clearly $D$ is also a facial cycle of % $G$. 
% Thus $B = {\cal  F}(G) - {D}$.                                
\ep
\\ 

Now we are ready to prove the following refinement of 
{\bf \ref{MacLane}}.
\bs 
%{\em (MacLane \cite{ML})}
\label{McLane} Let $G$ be a 2--connected graph.
\\[1ex]
$(a)$The following are equivalent:
\\[1ex]
\indent
$(a1)$ $G$ is planar and
\\[1ex]
\indent
$(a2)$ $G$ has a simple cycle basis.
\\[1ex]
$(b)$
Moreover,  if $G$ is not a cycle and $S$ is a simple basis of ${\cal CS}(G)$ then there exists an embedding $G_{\epsilon }$ of $G$ such that 
$S =  {\cal F}(G_{\epsilon }) \setminus \{F\}$ for some 
$F \in  {\cal F}(G_{\epsilon })$.
\es

\bp By {\bf \ref{Gplane}}, $(a1) \Rightarrow (a2)$.
We prove $(b)$ and $(a2) \Rightarrow (a1)$.
If $G$ is  a cycle, our claim is obviously true.
If $G$ is not a cycle, then our claim follows from 
{\bf \ref{2con}}.
\ep
\\

% Obviously {\bf \ref{MacLane}} is a Corollary of {\bf \ref{McLane}}.

\end{document}